\documentstyle[12pt]{article}
\newtheorem{fed}{Definition}[section]
\newtheorem{teo}[fed]{Theorem}
\newtheorem{lem}[fed]{Lemma}
\newtheorem{cor}[fed]{Corollary}
\newtheorem{pro}[fed]{Proposition}
\newtheorem{rem}[fed]{Remark}
\newtheorem{exa}[fed]{Example}

\def\dem{{\it Proof.\ }\rm}

\def\zC{\hbox{\rm C\hskip -5.4pt\vrule
height 8.0pt width 0.4pt depth -0pt\hskip 4.5pt}}
\def\zR{\hskip 2pt\hbox{\rm R\hskip -10.4pt{\rm I}
\hskip 2.5pt}}

\def\bull{\vrule height 1.0ex width 0.4ex depth -.1ex }
\def\inc{\subseteq}
\def\QED{\bull}
\def\inv{^{-1}}
\def\*A{\#\sb A}
\def\pa{P}

\def\H{{\cal H}}
\def\com{$(A, \cS)$ is compatible}

\def\Gpo{GL(\H)^+}
\def\ca{L(\H ) }

\def\cH{{\cal H}}
\def\glh{GL(\H)}
\def\cP{{\cal P}}

\def\cQ{{\cal Q}}
\def\cS{{\cal S}}
\def\cT{{\cal T}}
\def\cM{{\cal M}}
\def\cN{{\cal N}}
\def\PA{P_{A, \cS}}
\def\PB{P_{B, \cS}}
\def\QA{Q_{A,\cS} }
\def\csha{\Sigma (P, A)}
\def\PAS{\cP(A, \cS)}
\def\PBS{\cP(B, \cS)}
\def\rai{^{1/2}}
\def\mrai{^{-1/2}}
\def\api{\langle}
\def\cpi{\rangle}

\def\noi{\noindent}
\def\qa{\cQ}
\def\csta{C$^*$-algebra}

\def\bm{\left(\begin{array}}
\def\em{\end{array}\right)}
\def\ben{\begin{enumerate}}
\def\een{\end{enumerate}}
\def\beq{\begin {equation}}
\def\eeq{\end {equation}}
\def\barr{\begin{array}}
\def\earr{\end{array}}
\def\inv{^{-1}}
\def\pa{\cP}

\def\H{\cH}
\def\la{\lambda}

\date{}

\begin{document}
\title{Oblique projections  and Schur complements }
\author {G. Corach
\footnote{Partially supported by Fundaci\'on Antorchas,
CONICET (PIP 4463/96) , Universidad de Buenos Aires (UBACYT TX92 and TW49) and
ANPCYT (PICT97-2259)}\ ,  A. Maestripieri and D. Stojanoff
\footnote{Partially supported CONICET (PIP 4463/96),
Universidad de Buenos Aires (UBACYT TW49) and ANPCYT (PICT97-2259)}}

\maketitle

\noindent{Gustavo Corach, Depto. de Matem\'atica, FCEN-UBA, Buenos Aires,
Argentina}

\noindent{e-mail: gcorach@mate.dm.uba.ar}

\noindent{Alejandra Maestripieri, Instituto de
Ciencias, UNGS,  San Miguel, Argentina}

\noindent{e-mail : amaestri@ungs.edu.ar}

\noindent{Demetrio Stojanoff, Depto. de Matem\'atica, FCE-UNLP,  La Plata,  Argentina }

\noindent{e-mail: demetrio@mate.dm.uba.ar}
\bigskip


\begin{abstract}{
Let ${\cal H}$ be a Hilbert space, $L({\cal H} )$ the algebra 
of all bounded linear operators on ${\cal H}$ and 
$ \langle  ,  \rangle_A : {\cal H} \times {\cal H} \to \zC$ 
the bounded sesquilinear form induced by a  selfadjoint $A\in L({\cal H} ) $, 
$ \langle \xi, \eta \rangle_A = \langle A \xi , \eta \rangle \ , \quad \xi , 
\ \eta \in {\cal H} .$ 
Given $T\in L({\cal H} )$,  $T$ is $A$-selfadjoint if 
$AT = T^*A$. If ${\cal S}  \subseteq {\cal H}$ is a closed
subspace, we study the set of $A$-selfadjoint projections onto 
${\cal S} $,  
$$
\cP(A, {\cal S} )  = \{Q \in \ca: Q^2 = Q\ , \ \  
R(Q) = {\cal S}  \ , \ \ AQ = Q^*A\}
$$ 
for different choices of $A$, mainly under the hypothesis that $A\ge 0$.
There is a closed relationship between the $A$-selfadjoint projections 
onto ${\cal S} $ and the shorted operator (also called Schur complement) of $A$ 
to ${\cal S} ^\perp$. Using this relation we find several conditions which are
equivalent to the fact that $\cP(A, {\cal S} )  \neq \emptyset$, in particular 
in the case of $A\ge 0$ with $A$ injective or with $R(A)$ closed. 
If $A$ is itself a projection, we relate the set $\cP(A, {\cal S} ) $ with the existence 
of a projection with fixed kernel and range and we determine its norm.
}\end{abstract}

\newdimen\normalbaselineskip
\normalbaselineskip=16pt
\normalbaselines
\vskip 1truecm

\section{Introduction}

If $\H$ is a Hilbert space with scalar product $\api \ \cdot , \cdot \rangle$
and $L(\H)$ is the algebra of 
all bounded linear operators on $\H$, consider the subset 
$\cQ$ of $L(\H )$ consisting of all projections onto (closed) 
subspaces of $\H$ and the subset $\cP$ of $\cQ$ of all orthogonal 
(i.e., selfadjoint) projections. Every $Q \in \cQ \setminus \cP$ is
called an oblique projection. The structure of $\cQ$ and $\cP$ 
has been widely studied since the begining of the
spectral theory. 
In recent times, applications of oblique 
projections to complex geometry \cite{[PW1]}, statistics   \cite{[W1]}, 
\cite{[W2]} and wavelet theory \cite {[Al1]}, \cite {[Al2]}, \cite{[T]},
\cite{[T2]} have renewed the 
interest on the subject. The reader is also referred 
to \cite {[Bu]}, \cite{[HN]}. 

In \cite {[PW1]}, \cite{[ACS1]} there is an analytic study of the 
map which assigns to any positive invertible operator $A \in L(\H )$ and
any subspace $\cS$ of $\H$ the unique projection onto $\cS$ 
which is selfadjoint for the scalar product $\api \cdot , \cdot \cpi_A$ 
on $\H$ defined by $\api \xi , \eta \cpi_A = \api A\xi , \eta \cpi$ 
$(\xi ,\eta \in \H)$. 
In this paper we study the existence of projections onto $\cS$ 
which are selfadjoint for $\api \cdot , \cdot \cpi_A$ if $A$ is not necesarily
invertible. More precisely, if $\cS$ is a closed
subspace of $\H$ and $B : \H \times  \H \to \zC$ is a 
Hermitian sesquilinear form, consider the subsets of $\cQ$, 
$$ 
\cQ _{\cS} = \{ Q \in \cQ :  Q(\H ) = \cS  \}  \quad 
\hbox { (projections with range S) }
$$ 
and
$$
\cQ^B =  \Big{\{} Q \in \cQ :B(\xi , Q \eta ) = B( Q \xi , \eta ) \ , 
\ \hbox { for all } \ \xi , \eta \in \H \Big{\}} 
\quad 
\hbox { (B-symmetric projections).}
$$
The main theme of the paper is the 
characterization of the intersection of 
$\cQ _{\cS}$ and $\cQ^B$. We shall limit our study
to the case in which $B$ is bounded, so that, by Riesz' theorem, 
there exists a unique selfadjoint operator $A \in L(\H)$
such that 
$$
B(\xi,\eta)=B_A(\xi,\eta)=\api A\xi ,\eta \cpi ; 
$$
we search to characterize the set 
$$ 
\PAS = \cQ _{\cS} \cap \cQ^{B_A} .
$$
Observe that $\PAS$ has a unique element if $A$ is a positive 
invertible operator, but in general it can have 0, 1 or
infinite elements.  Even if we get a characterization of 
$\PAS$ in general, much more 
satisfactory results  can be obtained for a positive $A$ ($A\ge 0$, i.e.
$\api A\xi ,\xi \cpi \ge 0 $ for all $\xi \in \H$).
In this paper, a pair $(A, \cS )$ consisting of a bounded selfadjoint operator
$A$ and a closed subspace $\cS \inc \H$ is said to be $compatible$ if $\PAS$ 
is not empty.

The contents of the paper are the following: 

In section 2 we collect 
several known results we shall use later. We show  in this
section that if $Q \in \cQ$, $A \in L(H)$ and $R(QA) \inc R(A)$. Then
the unique operator $D\in L(\H )$ verifying that
$$
QA= AD \ , \quad \ker D = \ker QA \quad \quad and \quad \quad R(D) \inc
\overline{R(A^*)},
$$
(called the reduced solution of  $AX=QA$) satisfies also that $D^2
= D$, i.e., $D \in \cQ$.

In section 3, some characterizations of the compatibility 
of $(A, \cS)$ are given; some of them hold for general selfadjoint 
operators $A$, and others hold only 
for positive operators $A$.
Among other properties, it is shown that
an oblique projection $Q$ is $A$- seladjoint 
(if $A\ge 0$) if and only if $0\le Q^*AQ \le A$ (see Lemma \ref{RAR}).
We establish,  also for $A\ge 0$, that $\PAS $ is an affine manifold 
and we give a parametrization for it. When $(A, \cS)$ is compatible, a distinguished
element $\PA \in \PAS$ can be defined. It is shown that the norm of $\PA$ is
minimal in $\PAS$ (see Theorem \ref{losPA}).

In section 4 we consider the relationship between the compatibility of $(A, \cS)$ and 
some properties of the Schur complement.
M. G. Krein \cite{[K]} and W. N. Anderson and G. E. Trapp \cite{[AT]},  extended  
the notion of Schur complement of matrices to Hilbert space operators, 
defining what it is  called the
{\it shorted operator}.  We recall the
definition: if $ A \in \ca ^+ $, $\cS \inc \H$ is a closed 
subspace and $P = P_\cS$ is the
orthogonal projection onto $\cS$, then the set 
$$
\{ X \in \ca ^+ : X \le A \ \  \hbox { and } 
\ \  R(X)\inc \cS ^\perp  \}
$$ 
has a maximum (for the natural order relation in $\ca^+$), 
which is called the {\it shorted operator} of
$A$ to $\cS^\perp$. We shall denote it by $\Sigma (P, A)$. It is shown 
that, for any $Q \in \PAS$, 
the Schur complement $\Sigma (P, A)$ verifies that 
$$ 
\Sigma (P, A) = A(1-Q) 
$$
(see Proposition \ref{csqa}). We also show that 
\com \  if and only if, in the characterization
$$
\Sigma (P, A) = \inf \{ R^* A R \ : \ R\in \cQ ,\  \ker R = \cS \ \} ,
$$ 
due by Anderson and Trapp \cite{[AT]}, 
the infimum is, indeed, a minimum (see Corollary \ref{RAR2}).

In section 5 we consider the case of positive operators $A$
which are $injective$. Using properties of the 
shorted operator $\Sigma (P, A)$, new conditions equivalent to the 
fact that the pair \com \ are found. For example (see Proposition \ref{equiva}), 
it is shown that
$$
(A, \cS ) \ \ \hbox{ is compatible } \ \  \iff \cS ^\perp \inc R(A+\lambda (1-P)) , \ \ 
\mbox{ for some } \ \ \lambda >0 .
$$
In section 6 we consider the case of positive operators 
$A$ with $closed$ $range$. Among other equivalences, it is shown that 
\com \ if and only if  $\cS + \ker A$ is closed 
(see Theorem \ref{cerrado}).
As a consequence it is shown that all manifolds $\PBS$ for 
$R(B) = R(A)$ are "parallel" (see Corollary \ref{rango}).
So, in this sense, it suffices to study the case of the orthogonal
projection $Q = P_{R(A)}$. This case is studied in section 7, where we
show a formula for the norm of the projection $P_{Q,P} := P_{Q, \cS}$
in ${\cal P} (Q, \cS)$.
For example (see Proposition \ref{norma}), if $\ker Q \cap R(P) = \{0\}$, then 
$PQP \in GL(\cS) $ and 
$$
\|P_{Q,P} \|^2 = \|(PQP)\inv \| =(1-\|(1-Q)P\|^2)\inv .
$$
In case that
$R(P) \cap \ker Q = \{ 0 \} = R(Q) \cap \ker P$ (e.g., if $P$ and $Q$ are  in
{\it position p} \cite {[Di], [Da]} or {\it generic position } \cite{[Ha]}), 
$P_{Q,P}$ is the oblique projection given by
$$ 
\ker P_{Q,P} = \ker Q \quad \hbox{ and } \quad R(P_{Q,P}) = R(P) .
$$

\section{Preliminaries}
In this paper $\H$ denotes a Hilbert space, $\ca$ is the algebra of
all linear bounded operators on $\H$, $\ca^+$ is the subset
of $\ca$ of all (selfadjoint) positive operators, $GL(\H)$ is the group of all
invertible operators in $\ca$ and $GL(\H)^+ = \glh \cap \ca^+$ 
(positive invertible operators). For every $C \in
\ca$ its range is denoted by $R(C)$.

Denote by $\qa$ (resp. $\pa$) the set of all
projections (resp. selfadjoint projections) in $\ca$:
$$ \qa = \qa(\ca )= \{ Q \in \ca : Q^2 = Q\}  
\quad , \quad 
\pa = \pa(\ca ) = \{ P \in \qa : P = P ^* \}  .$$
The nonselfadjoint elements of $\qa$ will be called {\it oblique
projections}. 

Along this note we use the fact that every $ P \in \pa$ induces a
representation of elements of $\ca$ by $2 \times 2$ matrices. Under 
this representation
$P$ can be identified with
$$ \left(
\begin{array}{cc} I\sb {P (\H)}&0 \\ 0&0 \end{array} \right) =
\left( \begin{array}{cc} 1 &0 \\ 0&0 \end{array} \right)  $$
and all idempotents $Q$ with the same range as $P$ have the form
$$ Q= \left( \begin{array}{cc} 1&x \\ 0&0 \end{array} \right) $$
for some $x \in L(\ker P , R(P))$. 

Now we state the well known criterium due to Douglas \cite{[Do]} (see also
Fillmore-Williams \cite{[FW]}) about ranges and factorizations of operators:
\begin{teo}\label{douglas}
Let $A, B \in \ca$. Then the following conditions are equivalent:
\ben
\item $R(B) \inc R(A)$.
\item There exists a positive number $\lambda $ such that $BB^* \le \lambda AA^*$.
\item There exists $D \in \ca$ such that $B = AD$.

\een
Moreover, the operator $D$ is unique if it satisfies the conditions
$$
B= AD \ , \quad \ker D = \ker B \quad \quad and \quad \quad R(D) \inc
\overline{R(A^*)}.
$$
In this case $\|D\|^2 = \inf \{\lambda : BB^* \le \lambda AA^* \}$
and $A$ is called the {\bf reduced} solution of the equation $AX=B$.

\end{teo}
\begin{cor}\label{cuadrado}
Suppose that $Q \in \cQ$, $A \in L(H)$ and $R(QA) \inc R(A)$. Then
the reduced solution  $D\in L(\cH )$ of  $AX=QA$ satisfies that $D^2
= D$, i.e., $D \in \cQ$.
\end{cor}
\dem Note that $AD^2= QAD = Q^2A = QA$. Also 
$$ 
\ker QA = \ker D \inc \ker D^2 \inc \ker AD^2 =  \ker QA 
$$ 
and $R(D^2) \inc R(D) \inc \overline{ R(A^*)}$. 
Thus, $D^2$ is a reduced solution of $AX= QA$ and, by uniqueness, 
it must be $D^2 = D$, i.e. $D \in \cQ$ \quad
\QED

\section {$A$-selfadjoint projections, generic properties}

Throughout,  $\cS$ is a closed subspace of $\H$ and $P$ is the orthogonal
projection onto $\cS$. As we said in the introduction, we 
consider a bounded sesquilinear form $B = B_A: \H \times \H \to \zC$
determined  by a Hermitian operator $A \in \ca$:
$$
B_A(\xi, \eta ) = \api A \xi , \eta \cpi \ , \quad \xi , \ \eta \in \H  .
$$
This form induces the notion of $A$-orthogonality. For example, easy
computations show that the  $A$-orthogonal of $\cS$ is
$$
\cS^{\perp _A} :=
\{\xi: \ \api A \xi, \eta \cpi \  = 0 \ \ \forall \eta \in \cS\ \} =
A\inv (\cS ^\perp ).
$$
Given $T \in \ca$, an operator $W \in \ca$ is called an {\it A-adjoint}
of $T$ if
$$
B_A(T \xi, \eta)= B_A(\xi, W \eta), \quad \xi,
\ \eta \in \H,
$$
or, which is the same, if
$$
T^*A = AW .
$$
Observe that $T$ may have no $A$-adjoint, only one
or many of them. We shall not deal in this paper with the general
problem of existence and uniqueness of $A$-adjoint operators.
Instead, we shall study the existence and uniqueness of
$A$-selfadjoint projections, i.e., $Q \in \cQ$ such that $AQ = Q^*A$.
Among them, we are interested in those whose range is exactly $\cS$.
Thus, the main goal of the paper is the study of the set $$\PAS = \{Q
\in \cQ: R(Q) = \cS, AQ = Q^*A\}$$ for different choices of $A$. 
\begin{fed} \rm 
Let $A =A^* \in \ca $ and $\cS \inc \H$ a closed subspace.
The pair $(A, \cS )$ is said to be $compatible$ 
if there exists an $A$-selfadjoint projection 
with range $\cS$, i.e. if $\PAS $ is not empty.
\end{fed}
For general results on $A$-selfadjoint operators 
the reader is referred to the papers by 
P. Lax \cite{[L]} and J. Dieudonn\'e \cite{[D]};  
a recent paper by S.Hassi and K.  Nordstr\"om \cite{[HN]} contains 
many interesting results on $A$-selfadjoint projections. 
Some of the results of this section overlapp 
with  their work, but we include them because the methods 
used in our proofs are useful 
for the study of the case of a positive $A$, which is our main concern.

\begin{lem}\label{RAR}  Let $A =A^* \in \ca $ and $Q\in \cQ $.
Then  the following conditions are equivalent:
\ben
\item $Q$ verifies that $AQ =Q^*A$, i.e. $Q$ is $A$-selfadjoint.
\item $\ker Q \inc A\inv (R(Q)^\perp ) = R(Q)^{\perp_A}$.
\een
If $A\in \ca^+$, they are equivalent to 
\ben
\item [3.] $Q^*AQ \le A$.
\een
\end{lem}
\dem

\noi
{\bf 1 $\leftrightarrow$ 2: }
If $Q \in \PAS$ and $\xi, \eta \in \cH$, then
\beq\label{kerQ}
\api A\eta ,
Q \xi \cpi =\api Q^* A \eta ,  \xi \cpi = \api AQ \eta , \xi \cpi = \api Q \eta , A
\xi \cpi ,
\eeq
so $\ker Q  \inc A\inv (S^\perp )$. The converse can be proved in
a similar way.

\noi
{\bf 1 $\leftrightarrow$ 3: }
Suppose that $0\le Q^*AQ \le A$.
Then, by Theorem \ref{douglas}, the reduced solution $D$
of the equation  $A\rai X = Q^*A\rai $ satisfies $\|D\|\le 1 $ and, by
Corollary \ref{cuadrado}, $D^2 = D$. Thus, it must be $D^* = D$. Since
$Q^*A = A\rai DA\rai$, we conclude that $Q^*A = AQ$.
Conversely,  note that $AQ  = Q^* AQ \ge 0$ and, if $E= 1-Q$,
$AE = E^*AE$. Then $AQ  \le A$, because, for $\xi  \in \cH$, 
$$
\barr{rl}
\api AQ  \xi , \xi  \cpi & = \api AQ  \xi , Q  \xi  \cpi \\&\\
                    & \le \api AQ  \xi , Q \xi  \cpi + \api AE  \xi , E  \xi  \cpi \\&\\
                    & = \api AQ  \xi , \xi  \cpi + \api AE  \xi ,   \xi  \cpi = 
					\api A( Q  + E ) \xi , \xi  \cpi \\&\\
                    & = \api A \xi , \xi  \cpi  \quad \QED
\earr
$$

Throughout, we use the matrix representation determined by $P$.

\begin{pro}\label{PA} Given $A =A^* \in \ca$, the following conditions are equivalent:
\ben
\item  The pair \com \ (i.e. $\PAS$ is not empty).
\item $R(PA) = R(PAP)$.
\item If $A = \bm {cc} a & b \\ b^* & c \em $ then $R(b) \inc R(a)$.
\item $\cS + A\inv (\cS ^\perp ) = \H$.
\een
\end{pro}
\dem Note that $$ PA =  \bm {cc} a & b \\ 0 & 0 \em  \quad \hbox{
and } \quad PAP =  \bm {cc} a & 0 \\ 0 & 0 \em , $$ so $ R(a) =
R(PAP) \inc R(PA) = R(a) + R(b)$ and items 2 and 3 are equivalent.
On the other hand, for any  $Q \in \cQ$ it holds $R(Q) = \cS$ if
and only if $$ Q =  \bm {cc} 1 & x \\ 0 & 0 \em . $$ Easy
computations show that $Q^*A = AQ$  if and only if $ ax = b$, so
items 1 and 3 are equivalent by Theorem \ref{douglas}. Finally,
if $Q \in \PAS$ then, by Lemma \ref{RAR},
$\ker Q \inc A\inv (\cS ^\perp )$, which implies 4.
Conversely, if $\cS + A\inv (\cS ^\perp ) = \H$, and if $\cN$ is defined by
$\cN = \cS \cap A\inv (\cS ^\perp )$, then $\cS \oplus (A\inv (\cS
^\perp ) \ominus \cN) = \H$.  The projection $Q$ defined by
this decomposition of $\H$ verifies, again by Lemma \ref{RAR},
that $Q ^* A = A Q $ \quad \QED

\begin {rem}\rm \ben 
\item
As mentioned before, there exist operators $T\in \ca$ which do not
admit $A$-adjoint. In fact, the existence of an $A$-adjoint $W$ of
$T$ is equivalent to the existence of a solution of the equation 
$AW = T^*A$ and this is equivalent to $R(T^*A) \inc R(A)$. 
If $Q\in \cQ$, then the existence
of an $A$-adjoint of $Q$ is also equivalent to $R(A) = R(A) \cap \ker
Q^\perp + R(A) \cap R(Q)^\perp$.
\item We conjecture that the
existence of some $Q \in \cQ_\cS$ which admits $A$-adjoint is
equivalent to the fact that \com .
\een
\end{rem}

\begin{fed} \label{1d}
Let $A = A^* \in \ca$ and suppose that the pair \com .
If $A = \bm {cc} a & b \\ b^* & c \em $ and $d\in L(\cS^\perp , \cS)$
is the reduced solution of the equation $ax=b$, we define
the following oblique projection onto $\cS$:
$$ \PA :=  \bm {cc} 1 & d \\ 0 & 0 \em $$
\end{fed}

\begin{teo}\label{losPA}
Let $A = A^* \in \ca$  and suppose that \com .
Then the following properties hold:
\ben
\item $ \PA \in \PAS $.
\item $\PAS$ has a unique element (namely, $\PA$) if and only if
$\cS \oplus A\inv ( \cS^\perp) = \cH$.
\een
If $A\in \ca^+$, then
\ben
\item [3.] $A\inv(S^\perp)\cap \cS = \ker A \cap \cS := \cN$
\item[4.] $\PAS $ is an affine manifold and it can be parametrized as
$$\PAS = \PA \ + \ L(\cS^\perp, \cN),
$$
where $L(\cS^\perp, \cN)$ is viewed as a subspace of $\ca$. 
A matrix representation of this parametrization is
\beq\label{1zd}
\PAS \ni Q = \PA + z =
\bm {ccc} 1 &0& d \\ 0 & 1& z \\ 0&0 & 0 \em \barr{l} \cS\ominus \cN
\\ \cN \\ \cS ^\perp \earr
\eeq
with the notations of Definition \ref{1d}.
\item [5.] $\PA$ has minimal norm in $\PAS$:
$$ \|\PA \| = \min \{\ \|Q\| : Q \in \PAS \}. $$
Nevertheless, $\PA$ is not in general
the unique $Q\in \PAS$ that realizes the minimum norm.
\een
\end{teo}
\dem
\ben
\item Use the same argument as in the proof of Proposition \ref{PA}.
\item By Lemma (\ref{RAR}), if
$Q \in \cQ$ and $R(Q)= \cS$, then $Q\in \PAS $ if and only if
$\ker Q \inc A\inv (S^\perp )$. This clearly implies 2.
\item  
If $\xi  \in A\inv(S^\perp)\cap \cS$, then 
$\|A\rai \xi\| = \api A\xi , \xi  \cpi = 0$, so that $A\xi = 0$.
\item We have to show that
every element $Q\in \PAS$ can be written in an unique form as
$$
Q = \PA + z \ , \quad \hbox{ with } \ \ z\in L( \cS^\perp , \cN ) .
$$
If $A = \bm {cc} a & b \\ b^* & c \em
$, $Q = \bm {cc} 1 & y  \\ 0 & 0 \em $ with $y  \in L(\cS^\perp ,
\cS)$ and $d\in L(\cS^\perp ,
\cS)$ is the reduced solution  of the equation $ax=b$,
then $Q\in \PAS$ if and only if  $ay = b$ if and only if  $a(y-d)
= 0$. Therefore, if $z = y-d \in L(\cS^\perp ,
\cS)$, then $Q\in \PAS$ if and
only if $Q= \PA + z$ and $R(z) \inc \ker a $.  But
$$
\ker a = \cS \cap \ker PAP = \cS \cap \ker A = \cN.
$$
Concerning the matrix representation, note that, by Theorem \ref{douglas},
$$R(d) \inc \overline {R(a)} = (\ker a) ^\perp = \cS \ominus \cN .
$$
\item
If $Q \in \PAS$ has the matrix form given in equation (\ref{1zd}), then
$$
\|Q\|^2   =
1+ \left\| \bm {ccc} 0 &0& d \\ 0 & 0 & z \\ 0&0 & 0 \em \right\|^2
\ge  1+\|d \|^2 = \|\PA\|^2 .
$$
Choose $d\in L(\cS^\perp , \cS)$ such that $\|d\|= 1$, $R(d) = 
\overline {R(d)} \neq \cS$ and $\ker d \neq \{0\}$. Then the matrix
$$
A = \bm{cc} P_{R(d)} & d \\ d^* & 1 \em \ge 0 , 
$$
$\cN = \ker A \cap \cS = \cS \ominus R(d)$ and $d$ is the reduced solution
of $P_{R(d)}x = d$.
Let $z \in L(\ker d , \cN )$ with $0 < \|z\|\le 1$; then
the projection $Q= \PA + z$ as in  equation (\ref{1zd}) satisfies $Q\in \PAS$, 
$\|Q\|=\|\PA\|= \sqrt {2}$ and $Q \neq \PA$ \quad \QED
\een

\section{Schur complements and $A$-selfadjoint projections}

As before, let $P \in \cP$ be the orthogonal projection onto the closed subspace
$ \cS \inc \cH$.
Every $A \in \Gpo$ defines a scalar product on $\H$ which is
equivalent to $\api , \cpi $, namely
$$
\api \xi, \eta \cpi_A = \api A\xi, \eta \cpi , \quad
\xi, \eta \in \H .
$$
The unique projection $\PA$ onto $\cS$ which is
$A$-orthogonal, i.e., $A$-selfadjoint, is uniquely determined by
$$
\PA = P(1 + P - A^{-1}PA)^{-1} = P(PAP + (1-P)A(1-P))^{-1}A.
$$
Observe that $\PA = A^{-1}\PA^*A$, because $A$ is invertible.  In
particular, in this case the set $\PAS$ is a singleton. Analogously,
there exists a unique projection $\QA$ which is $A$-orthogonal and
has kernel $\cS$: $\QA = 1 - \PA$. Notice that $A\QA = \QA^*A$.

Consider the map 
$$ 
\Sigma : \cP \times \Gpo \to \ca ^+ \ , \quad \hbox{ defined by }
\quad \Sigma (P, A)= A\QA = \QA^*A .
$$
If $A \in \Gpo$ has matrix representation
$ A = \bm {cc} a & b \\ b^* & c \em  ,$ then
$$\PA = \bm {cc} 1 & a^{-1}b \\ 0 & 0 \em ,
\quad  \QA = \bm {cc} 0 & -a^{-1}b \\ 0 & 1 \em , \
\hbox{ and } \
\Sigma(P, A) = \bm {cc} 0 & 0
\\ 0 & c  - b^*a\inv b \em .
$$
This reminds us the {\it Schur complement}. Recall that, given a square matrix
$M= \bm {cc} a & b \\ c & d \em $, with $a$ and $d$ square blocks, a 
Schur complement of $a$ in $M$ is $d-ca'b$, where $a'$ is a generalized inverse 
of $a$. The reader is referred to \cite{[Co]} and \cite{[Ca]} for 
concise surveys on the subject. This notion has been extended to 
positive Hilbert space operators by M. G. Krein \cite{[K]} and, 
later and independently, by  W. N. Anderson and G. E. Trapp
\cite{[AT]} defining what is called the {\it shorted operator}:
if $ A \in \ca ^+ $ then the set 
$$
\{ X \in \ca ^+ : X \le A \quad \hbox { and } \quad R(X)\inc \cS ^\perp  \}
$$
has a maximum (for the natural order
relation in $\ca^+$), which is called the {\it shorted operator} of
$A$ to $\cS^\perp$.

Next we collect  some results of Anderson-Trapp and E. L. Pekarev 
\cite{[PEKA]} which are relevant in this paper. 
Observe that the first item allows us to extend the 
map $\Sigma$ to $\cP \times \ca^+.$

\begin {teo} \label{csqa0}
Let $A \in \ca ^+ $ with matrix representation $$A = \bm {cc} a & b
\\ b^*  & c \em.$$
\ben
\item If $A$ is invertible, then $\Sigma (P, A)$ coincides with the
shorted operator of $A$ to $\cS^\perp$. We shall keep the notation
$\Sigma (P, A)$ for the shorted operator of $A$ to $R(P)^\perp$ for
every pair $(P, A) \in \cP \times \ca^+$.
\item $R(b) \inc R(a\rai)$ and if $d \in \ca$ is the reduced solution
of the equation $a\rai \ x = b$ then $$ \csha = \bm {cc} 0 & 0
\\ 0 & c  - d^*d \em
$$
\item If $\cM = A\mrai (\cS^ \perp )$ and $P_\cM$ is the orthogonal projection
onto $\cM$ then $$\csha = A\rai P_\cM A\rai . $$
\item $\csha$ is the infimum of the set $ \{ R^* A R \ : \ R\in \cQ ,\  \ker R = \cS \
\}$; in general, the infimum is not attained.
\item $R(A) \cap \cS^\perp \inc R(\csha )\inc R(\csha \rai ) = R(A\rai ) \cap
\cS^\perp$; in general, the inclusions are strict.
\een
\end{teo}
The reader is referred to \cite{[AT]} and \cite{[PEKA]} for proofs of these facts.
We prove now that the infimum of item 4 is
attained if and only if \com \ by relating the notions of
shorted operators and  $A$-selfadjoint projections (when there is one). 
As a consequence, we complete item 5 of the last theorem in case that 
\com .

\begin {pro} \label{csqa}
Let $A \in \ca ^+ $ such that the pair \com . Let
$E\in \PAS$ and $Q =1-E $. Then
\ben
\item $\csha = AQ  = Q^* AQ $.
\item $ \csha  = \min \{ R^* A R : R\in \cQ , \ker R = \cS \}$
and the minimum is attained at $Q$.
\item $R(\csha ) = R(A) \cap S^\perp $.
\een
\end{pro}
\dem
\ben
\item Note that $0\le AQ  = Q^* AQ \le A$, by Lemma \ref{RAR}.
Also $R(AQ ) = R(Q^*A) \inc R(Q^*) = \cS^\perp$.
Given $X \le A $ with $R(X) \inc \cS^\perp $, then, since
$\ker Q  = \cS$, we have that
$$
X = Q^* X Q  \le Q^* A Q
= AQ ,
$$
where the first equality can be easily checked because
$X$ has the form $  \bm {cc} 0 & 0 \\ 0 & x \em$.

\item By item 1, $Q^*AQ = \csha $ and $\ker Q = \cS$. So the
minimum is attained at $Q$ by Theorem \ref{csqa0}.
\item Clearly the equation $\csha = A Q $ implies that
$R(\csha ) \inc R(A) \cap S^\perp $. The other inclusion always holds by
Theorem \ref{csqa0} \quad \QED
\een

\begin{cor}\label{RAR2}If $A \in \ca ^+ $  the
following conditions are equivalent: \ben
\item  The pair \com .
\item The set $\{ S^* A S : \ S\in \cQ ,  \ \ker S =  \cS \} $
attains its minimum at some projection $R$.
\item There exists $R\in \cQ $ such that
$\ker R = \cS$ and $R^*AR \le A $.
\een
\end{cor}
\dem
\noi
{\bf 1 $\rightarrow$ 2: } Follows from Proposition \ref{csqa}. 

\noi
{\bf 2 $\rightarrow$ 3: } Follows from Theorem \ref{csqa0}. 

\noi
{\bf 3 $\rightarrow$ 1: } By Lemma \ref{RAR}, any projection $R$ such that
$R^*AR \le A$ verifies that $AR =R^*A$. If also $\ker R = \cS$, then $ 1-R
\in \PAS $ \QED

In the next sections we shall study the
existence of $A$-selfadjoint projections onto a closed subspace $\cS$,
under particular hypothesis on the positive operator $A$.

\section{$A$-selfadjoint projections: the injective case}
As before, let $P \in \cP$ be the orthogonal projection onto $ \cS $.
In this section we study the case of injective operators $A\in \ca^+ $.
We define the notion of $A$-admissibility for $\cS$, in terms of the 
shorted operator $\csha$. This notion is shown to be strictly weaker 
that compatibility for the pair $(A, \cS)$. Under the assumption of
$A$-admissibitity for $\cS$, the fact that \com \ becomes equivalent
to the equality $R(\csha) = \cS^\perp \cap R(A)$ (see 
item 5 of Theorem \ref{csqa0} and item 3 of Proposition \ref{csqa}).

\begin{lem}\label{admissible}
Given $A\in \ca ^+$ which is injective and $\cM = A\mrai (\cS
^\perp )$, the following conditions are equivalent: \ben
\item $\ker \csha = \cS$
\item $ \cM ^\perp \cap A\rai (\cS^\perp ) = \{ 0 \} $.
\item $\overline{\cS ^\perp \cap R(A\rai )} = \cS ^\perp $
\item  $\cS = ( \cS ^\perp \cap R( A\rai ) )^\perp $
\item $\cS = \cT ^\perp $ for some subspace $\cT \inc R(A\rai )$
\een
\end{lem}
\dem 
\noi
{\bf 1 $\rightarrow$ 2: }
Recall that $\csha = A\rai P_\cM A\rai $ 
by Theorem \ref{csqa0}. Using that $A\rai $ is injective, we deduce 2. 

\noi
{\bf 2 $\rightarrow$ 3: } Suppose that 3 is false. 
Let $ \xi  \in \cS ^\perp \ominus
\overline{\cS ^\perp \cap R(A\rai ) }$, $\xi  \neq 0$. 
If $\eta \in \cM $ then $A\rai \eta \in \cS ^\perp 
\cap R(A\rai )$ and 
$$ 
\api A\rai \xi , \eta \cpi = \api \xi , A\rai \eta \cpi = 0. 
$$ 
Thus, $A\rai \xi \in \cM ^\perp
\cap A\rai (\cS^\perp )$ and  $A\rai \xi \neq 0$, which contradicts 2.

\noi
{\bf 3 $\rightarrow$ 4 $\rightarrow$ 5: } It is clear.

\noi
{\bf 2 $\rightarrow$ 3: } 
If $\cS = \cT^\perp $ with $\cT \inc R(A\rai )$, 
then $\cS ^\perp = \overline {\cT}$ and $$ \cT \inc
R(A\rai ) \cap \cS^\perp = R(\csha \rai ). $$ So $R(\csha \rai )$ is
dense in $\cS^\perp$ and $\ker \csha = \ker \csha\rai = R(\csha \rai
)^\perp = \cS$ \quad \QED

\begin{fed}\rm We shall say that $\cS$ is $A$-{\it admissible} if
any of the conditions of Lemma \ref{admissible} is verified.
\end{fed}

\begin{lem}\label{1->2}
If $ A\in \ca ^+$is injective and \com , then $\cS$ is $A$-admissible.
\end{lem}
\dem Let $E \in \PAS$ and $Q = 1-E$. Then, by
Proposition \ref{csqa}, $\csha = A Q $ and $\ker \csha = \ker Q  = \cS$ \quad \QED

\begin{rem}\rm If \com \ then a condition which is
stronger than $A$-admissibili\-ty
is verified. Indeed, $\csha = A (1-\PA) $ implies that $\ker \csha =
\cS$. But in this case, 
$R(\csha )  \inc R(A) \cap \cS^\perp $ which must be dense in $\cS^\perp $.
Note that $R(A\rai )$  strictly contains $R(A)$ if $R(A)$ is not
closed.

Nevertheless, we restrict ourselves to the weaker notion of
$A$-admissibility because under the hypothesis that $\cS$ is 
$A$-admissible, the conditions \com \ and 
$R(\csha) \inc R(A)$ become equivalent. 
Observe that $R(\csha) \inc R(A)$  is false in general 
(recall item 5 of Theorem \ref{csqa0} and item 3 of Proposition \ref{csqa})
\end{rem}

\begin{pro}\label{equiva}
If $ A\in \ca ^+$ is injective then the following conditions are
equivalent: 
\ben
\item  The pair \com .
\item
\ben
\item [i)] $\ker \csha = \cS$ (i.e. $\cS$ is $A$-admissible) and
\item [ii)] $R(\csha) \inc R(A)$. 
\een
\item $\cS$ is $A$-admissible and,
if $\cM = A\mrai (\cS ^\perp )$, then $P_\cM A P_\cM \le \mu A$ for some
$\mu > 0$.
\item $\cS ^\perp \inc R(A+\lambda (1-P))$ for some (and then for any) $\la > 0$.
\een
\end{pro}
\dem

\noi
{\bf 1 $\to$ 2: } By Lemma \ref{1->2}, $\cS $ must be $A$-admissible. If
$\QA  = 1-\PA $, then $\csha = A \QA $ and 2 follows.

\noi
{\bf 2 $\to$ 3: } If $R( A\rai P_M A\rai ) \inc R(A)$ then
$R( P_M A\rai ) \inc R(A\rai )$ and 3 holds.

\noi {\bf 3 $\to$ 1: } Note that $P_M A P_M \le \mu A$ if and only if $R(P_M
A\rai ) \inc R(A\rai)$ if and only if there exists a unique $F \in \ca$ such
that $A\rai F = P_M A\rai $, $\ker (P_M A\rai ) \inc  \ker F $ and
$R(F) \inc R(A\rai)$. We shall see that $1-F \in \PAS $. Indeed, $F^2 = F$
by Corollary \ref{cuadrado}. $F$ is $A$-selfadjoint because $AF =
A\rai P_M A\rai = \csha$ which is selfadjoint. Finally,  $\ker F =
\cS$. Indeed, $AF = \csha$, so $\ker F = \ker \csha = \cS$ because
$\cS $ is $A$-admissible.

\noi {\bf 4 $\leftrightarrow$  1: } Using Proposition \ref{PA},
we know that the fact that \com \ 
only depends on the first row $PA$ of $A$. Therefore we can freely change $A$ by
$A+\lambda (1-P)$, for $\la > 0$.
In this case conditions 2 can be rewritten as condition 4,
since $\Sigma (P, A+\lambda (1-P) ) = \csha +\lambda (1-P) $.
\quad \QED

\begin{exa} \rm Given a positive injective operator $A \in L(\cH)$ with non-closed
range, it is easy to show that 
there exists $\xi  \in R(A\rai ) \setminus R(A)$. Let $P_\xi $ be the
orthogonal projection onto the subspace generated by $\xi $. Then $R(P_\xi ) 
\inc R(A\rai )$, so that, by Douglas' theorem, $P_\xi  \le \lambda A$
for some positive number $\lambda$ which we can suppose equal to
$1$, by changing $A$ by $\lambda A$. 
It is well known that this implies that the operator $B \in
L(\H \oplus \H)$ defined by $$ B = \bm {cc} A & P_\xi 
\\ P_\xi  & A \em  $$ is positive.  Let $\cS = \cH_1 = \H \oplus 0$. Then $\cS^\perp =
\cH_2 = 0 \oplus \H$. We shall see that 
$B$ is injective, $\cH_1$ is $B$-admissible, moreover
$\cH_2 \cap R(B)$ is dense in $\cH_2$, but $\PBS$ is empty.

Indeed, it is clear that $B$ does not verify condition 3 of Proposition \ref{PA},
so $\PBS$ is empty. Let $D$ be the reduced solution of $P_\xi  = A\rai X $.
Then $\Sigma (P, B) =
A- D^*D$. Note that $\ker D = \ker P_x $ implies $DP_\xi  = D$. So
$D^*D = P_\xi  D^*D$. Then, if $0\oplus \eta \in \ker \Sigma (P, B)$,
$$
A\eta =  D^*D \eta =  P_\xi D^*D\eta = \lambda \xi    \ \hbox{ for some } \ \la \in \zR
\quad \Rightarrow \quad \eta = 0
$$
because $\xi  \notin R(A)$ and $A$ is injective.
So $\ker \Sigma (P, B) = \cS$ and $\cH_1$ is $B$-admissible. Also
\beq\label{kerB}
B(\omega \oplus \eta) \in 0 \oplus \cH   \iff A\omega + P_\xi  \eta = 0 \iff \omega
=0 \ \hbox{ and } \ \eta \in \{\xi \}^\perp.
\eeq
Then $R(B) \cap \cH_2 = \{ B(0 \oplus \eta ): \eta \in \{\xi \}^\perp\} = 0 
\oplus A(\{\xi \}^\perp ) $.
We shall see that $A(\{\xi \}^\perp )$ is dense in $\cH$. Indeed, if $\zeta \in
[A(\{\xi \}^\perp )]^\perp$, then
$\api \eta , A\zeta \cpi = \api A\eta , \zeta \cpi =  0  $
for all $ \eta \in \{\xi \}^\perp $. So $A\zeta = \mu \xi  $ for some $\mu \in \zR$.
As before this implies that $  \zeta = 0$.
Finally, the injectivity of $B$ can be easily deduced from 
equation (\ref{kerB}).
\end{exa}

\section{$A$-selfadjoint projections: the closed range case}

As before we fix $P\in \cP$ with $R(P)=\cS$.
In this section $A$ denotes a positive operator  with closed range. 
We shall see that, in this case, the fact that 
\com \ depends only on the angle between $\ker A$ and $\cS$, i. e. 
\com \ if and only if $\ker A + \cS $ is closed. To establish the link between 
compatibility and the angle condition, we need to determine when $R(PAP)$ is closed.
This is done in the following Lemma:

\begin{lem}\label{PAP}
It holds that
$$
\overline{R(PAP)} = \cS \cap (\cS \cap \ker A)^\perp 
$$
and that $R(PAP)$ is closed if and only if the subspace $\ker A + \cS $ is
closed.
\end{lem}
\dem
First note that $\ker PAP = \{ x \in \cH : \api PAP x,x\cpi =0 \}
= \{  x \in \cH : Px \in \ker A \} = \ker AP .$
So $\ker PAP = \cS^\perp \perp (\cS \cap \ker A )$. Therefore
$$
\overline{R(PAP)} = (\ker PAP )^\perp = \cS \ominus (\cS \cap \ker A )
= \cS \cap (\cS \cap \ker A )^\perp  := \cM .
$$
Clearly $\cM \cap \ker A = \{ 0 \}$. Suppose that
$\cN = \ker A + \cS = \ker A \oplus  \cM $ is closed. Let $Q$
the projection from $\cN$ onto $\cM$ with $\ker Q = \ker A$; observe that $Q$
is bounded. If $Q=0$ then $\cM = \{0\}$, $\cS \inc \ker A$ and $PAP = 0$.
If $\cM \neq \{0\}$, given $\xi  \in \cM$, let $\eta \in R(A)$ such that $A\eta  = A\xi $
($A$ is invertible in $R(A)$). Clearly $\eta  = \xi + \zeta$ with $\zeta \in \ker A$.
Then $\eta  \in \cN$, $Q\eta  = \xi $ and $\|\xi  \| \le \|Q\| \ \|\eta \|$. Therefore
$$
\api PAP \xi , \xi \cpi = \api A\xi ,\xi  \cpi = \api A\eta ,\eta  \cpi
\ge \lambda \|\eta \|^2 \ge \lambda \|Q\|^{-2} \|\xi \|^2 .
$$
for some $\lambda > 0$, since $A|_{R(A)}$ is bounded from below.
Conversely, if $R(PAP)$ is closed then $R(PAP) = \cM$.
Then there exists $\mu >0$ such that
$\api A\xi ,\xi  \cpi = \|A\rai \xi \|^2 \ge \mu \|\xi \|^2 $ for $\xi  \in \cM$ and
$A\rai (\cM )$ is closed. So $\cN = A\mrai ( A\rai (\cM ))$
must be also closed \ \ \ \QED

\begin {teo}\label{cerrado}
If $A \in \ca ^+$ has closed range then the following conditions are
equivalent: \ben
\item  The pair \com .
\item $R(PAP)$ is closed.
\item $\cS + \ker A$ is closed.
\item $R(PA)$ is closed.
\item $\cS^\perp + R(A)$ is closed.
\item $R(AP)= A(\cS )$ is closed.

\een
\end{teo}
\dem
By Lemma \ref{PAP} conditions 2 and 3 are equivalent.

{\bf \noi 2 $\to $ 1:} Let $A = \bm {cc} a & b \\ b^* & c \em $
in terms of $P$. Note that $a = PAP$, so $R(a)$ is closed.
Therefore, since $A\ge 0$, $R(b) \inc R(a\rai ) = R(a)$.
Then  \com \  by Proposition \ref{PA}.

{\bf \noi 1 $ \to $ 3:} Suppose that  \com .
Let $\PA \in \PAS$ and let $\QA = 1-\PA$. Then
$$
\barr{rl}
\ker A & \inc \ \ker (\QA ^{^*} A ) = \ker (A \QA  ) \\&\\
       & = \{\xi  \in \cH : \QA  \ \xi  \in \ker A \} \\&\\
       & = \cS \oplus (\ker A \cap R(\QA ) ) \\ &\\
       & \inc \ker A +\cS .
       \earr
$$
Therefore $\ker A +\cS = ( \ker A \QA )$
which is closed.

{\bf \noi 4 $ \leftrightarrow $ 5:} This is an easy consequence of the identity
$$R(A) + \cS^\perp = P\inv [P(R(A))] = P\inv [R(PA)].
$$

{\bf \noi 3 $ \leftrightarrow $ 5:} In fact, it holds in general that
the sum of two closed subspaces is closed if and only if the sum of
their orthogonal complements is closed (see \cite{[De]}).

{\bf \noi 4 $ \leftrightarrow $ 6 :} It is a general fact that $R(C)$
is closed if and only if $R(C^*)$ is closed. \quad \QED

\begin{rem}\rm Conditions 3, 4 and 5, 6 are known to be equivalent, since
$R(P)= \cS$ and $\ker P = \cS^\perp$ (see Thm. 22 of \cite{[De]}).
They are also equivalent to, for example, the angle condition
$$
c(\cS, \ker A ) < 1,
$$
where $c(\cS, \cT )$ is the cosine of the
Friedrichs angle between  the two subspaces $\cS ,  \ \cT$,
defined by:
\beq\label{angle}
c(\cS, \cT ) = \sup \{ |\api \xi , \eta  \cpi | : \xi \in \cS\cap (\cS \cap \cT)^\perp, \
\|\xi \|\le 1, \ \eta \in \cT\cap (\cS \cap \cT)^\perp, \ \|\eta \|\le 1 \}
\eeq
Also Lemma \ref{PAP} can be deduced from the results of \cite{[De]}.
\end{rem}

\begin {cor}\label{rango}
For every $A \in \ca^+$ with closed range, the following conditions
are equivalent:
\ben
\item The pair \com .
\item For all $B \in \ca^+$ with $R(B) = R(A)$, 
the pair $(B, \cS )$ is compatible.
\item The pair $(P_{R(A)}, \cS)$ is compatible, if $P_{R(A)}$ denotes
the orthogonal projection onto the closed subspace $R(A)$
\een
Moreover, if $B \in \ca^+$ and $R(B) = R(A)$, then the affine manifolds
$\PAS $ and $\PBS$ are "parallel", i.e.
\beq\label{paralelas}
\PBS = ( \PB  - \PA) \  + \ \PAS.
\eeq
\end{cor}
\dem If $R(B)= R(A)$ then $\ker B = \ker A = \ker P_{R(A)}$ 
and, by Theorem \ref{cerrado}, the three conditions
are equivalent. Equality (\ref {paralelas}) follows from the 
parametrization given in Theorem \ref{losPA}, since
$$
A\inv (\cS^\perp ) \cap \cS = \ker A \cap \cS= \ker B \cap \cS =
B\inv (\cS^\perp ) \cap \cS \quad \QED
$$
Condition 3 is an invitation to consider the sets ${\cal P}(Q, \cS)$ for
$Q \in \cP$, which we study in the next section.

\section{The case of two projections}

In this section we shall study the case in which $A$ is an 
orthogonal projection, i.e., $A = Q \in \cP$. 
Then, by Theorem \ref{cerrado} (items 3 and 6),
$\ker Q + R(P)$ is closed if and only ${\cal P}(Q, R(P))$ is not empty.
In this case we shall denote by $P_{Q,P}$ the projection
$ P_{Q, R(P)}$ of Definition \ref{1d}.
In the following theorem we collect several conditions which
are equivalent to the existence of $P_{Q,P}$. Notice, however, that 
the equivalence of items 3 to 10 
can be deduced from results by R. Bouldin \cite {[Bo]}
and S. Izumino \cite {[Iz]}; a nice survey on this and related subjects 
can be found in  \cite {[De]}. Observe that Theorem \ref{cerrado} 
provides alternative proofs of some of the equivalences.

\begin{teo}\label{largo}
Let $P, Q \in \cP$ with $R(P) = \cS$ and $R(Q) = \cT$. The
following are equivalent:
\ben
\item $(Q, \cS) $ is compatible.
\item $(P, \cT) $ is compatible.
\item $\ker Q + R(P) $ is closed.
\item $\ker P + R(Q) $ is closed.
\item $R(PQ)$ is closed.
\item $R(QP)$ is closed.
\item $R(1-P+ Q)$ is closed.
\item $R(1-Q + P)$ is closed.
\item $ c(\cS , \cT^\perp) =c(\cT , \cS^\perp) <1$.
\een
If $\ker Q \cap \cS = \{0\}$, they are equivalent to 
\ben \item[10.] $ \|(1-Q)P \|<1 $.
\een
\end{teo}
\dem

{\bf \noi 1 $ \leftrightarrow $ 2 $ \leftrightarrow $ 3:} Follows from Theorem
\ref{cerrado}.

{\bf \noi 3 $ \leftrightarrow $ 4 $ \leftrightarrow $ 9:} Follows from theorem 13
of \cite{[De]}.

{\bf \noi 3 $ \leftrightarrow $ 6 and 4 $ \leftrightarrow $ 5:} Follows from theorem 22
of \cite{[De]}.

{\bf \noi 5 $ \leftrightarrow $ 7 and 6 $ \leftrightarrow $ 8:} Follows from  2.5
of \cite {[Iz]}

{\bf \noi 3 $ \leftrightarrow $ 10} Follows from theorem 13
of \cite{[De]} \quad  \QED

\bigskip
\noi
Suppose that any of the conditions of Theorem \ref{largo} is verified by $P, Q \in \cP$. 
As a final result, we shall compute $\|P_{Q,P} \|$. First, we assume that
$ \ker Q \cap R(P) = \{0\}$:

\begin{pro}\label{norma}
Let $P, Q \in \cP$. Denote $R(P) = \cS$. Suppose that
$\ker Q \cap \cS = \{0\}$ and $\ker Q + \cS $ is closed. Then
$Q|_\cS $ is invertible in $L(\cS , Q(\cS ))$,
$PQP $ is invertible in $L(\cS )$ and
$$
\|P_{Q,P} \| = \|(Q|_\cS )\inv \| = \|(PQP)\inv \|\rai = ( 1 - \|(1-Q)P\|^2 )\mrai
$$
\end{pro}
\dem
Using Theorem \ref{largo}, we know that $\|(1-Q)P\|<1$. Then
$$
\|P-PQP\| = \|P(1-Q)P\| = \|(1-Q)P\|^2 < 1 ,
$$
showing that $PQP$ is invertible in $L(\cS )$. 
On the other hand consider $  Q|_\cS  : \cS \to Q(\cS ) $. 
By   Theorem \ref{cerrado},
$Q(\cS ) $ is closed, so $Q|_\cS $ is invertible in $L(\cS, Q(\cS ) )$.

If $P_{Q, P} = \bm {cc}  1& d \\0 & 0 \em $, then $\|
P_{Q, P} \|^2 = 1+\|d\|^2$. Recall that $d$ is the reduced solution
of the equation $PQP X = PQ(1-P)$. So, by Theorem \ref{douglas},
$$
\barr{rl}
\|d\|^2 & = \inf \{ \la >0 : PQ(1-P)QP \le \la PQPQP\} \\&\\
        & = \inf \{ \la >0 : PQP \le (1+\la ) (PQP)^2 \}\\&\\
        & = \inf \{ \la >0 : P  \le (1+\la ) PQP \}
		 = \inf \{ \la >0 : (PQP)\inv  \le (1+\la ) P \}\\&\\
        & = \| (PQP)\inv \| -1 .        \earr
$$
So $\|P_{Q,P} \| ^2 =  \| (PQP)\inv \|$. Note also that
$$
P  \le (1+\la ) PQP \iff \|\xi \|^2 \le (1+\la ) \api PQP \xi , \xi \cpi = (1+\la ) \|Q \xi \|^2
\quad \hbox{ for all } \  \xi  \in \cS .
$$
Taking infimum over $\la$, we get $ \| P_{Q, P} \| = (1+\|d\|^2)\rai = 
\|(Q|_\cS )\inv \| . $

It is easy to see that, if $0< A \le I $ in $L(\cH)$, then
$ \|I-A \| = 1-\|A\inv \|\inv $. Applying this identity
to $PQP$ in $L(\cS)$ we get
$$
\|(PQP)\inv \| = (1 - \| P- PQP \|)\inv = (1 - \| P(1-Q)P \|)\inv
= ( 1 - \|(1-Q)P\|^2 )\inv \quad \QED
$$

\begin{rem} \label{remPQ}\rm
Let $P, Q \in \cP$ with $R(P) = \cS$ and $R(Q) = \cT$ and suppose that
any of the conditions of Theorem \ref{largo} hold.
By Proposition \ref{losPA}, 
$$\barr{rl}
\ker P_{Q,P}&
= Q\inv (\ker P) \ominus (\ker Q \cap R(P)) \\&\\
& = (\ker Q + R(Q)\cap \ker P)\ominus (\ker Q \cap R(P)). \earr
$$
Therefore, in the case that
\beq\label{generic}
R(Q)\cap \ker P = \{0\} =  \ker Q \cap R(P)
\eeq
(e.g., if $P$ and $Q$ are  in
{\it position p} \cite {[Di], [Da]}
or {\it generic position } \cite{[Ha]}) 
we can conclude that $P_{Q,P}$ is the projection given by
$$ \ker P_{Q,P} = \ker Q \quad \hbox{ and } \quad R(P_{Q,P}) = R(P)
$$
Then $\cS \oplus \ker Q =\cH$ and $P_{Q, P}$ is the oblique projection
given by this decomposition of $\cH$. In this case,  formula 
$\|P_{Q,P} \| = ( 1 - \|(1-Q)P\|^2 )\mrai$ has been proved by 
Ptak in \cite{[Pt]} (see also \cite{[Bu]}).
\end{rem}

\begin{teo} Let $P, Q \in \cP$ which verify that $\ker Q + R(P) $ is closed.
Denote by $\cN =  \ker Q \cap R(P) $, $\cM = R(P)\ominus \cN$ 
and $P_0 = P_\cM$. Then 
\ben
\item ${\cal P}(Q, \cM )$ has only  one element, namely $P_{Q, P_0}$,  
\item $P_{Q, P} = P_\cN + P_{Q, P_0}$ and 
\item $\|P_{Q,P} \| = \|P_{Q,P_0} \| = (1-\|(1-Q)P_0\|^2)\mrai$.
\een
\end{teo}
\dem If $\cN = \{ 0 \}$, we can use Proposition \ref{norma}. 
Assume now that $\cN$ is not trivial. Then, by the results of section 3,
we get the matrix form
$$
P_{Q, P} = \bm {ccc} 1 &0& 0 \\ 0 & 1& d \\ 0&0 & 0 \em \barr{l}
\cN \\ \cM \\ \ker P \earr   .
$$
Denote
\beq\label{TT}
T= P_{Q, P}- P_\cN = \bm {ccc} 0 &0& 0 \\ 0 & 1& d \\ 0&0 & 0 \em \barr{l}
\cN \\ \cM \\ \ker P \earr  .
\eeq
We must show that $T = P_{Q, P_0}$. Note that $\ker Q \cap \cM = \{ 0\}$, so 
${\cal P}(Q, \cM )$ has, at most, one element. On the other hand, $T^2 =T$ and $R(T) = \cM$
by equation (\ref{TT}). Also
$$
T^*Q = (T^*+P_\cN)Q= P_{Q, P}^*Q = QP_{Q, P} = Q(P_\cN +T) = QT ,
$$
because $QP_\cN = 0$. So, $T = P_{Q, P_0}$ as claimed. 
By equation (\ref{TT}) and Proposition \ref{norma},
$$
P_{Q, P} = P_\cN + P_{Q, P_0} \quad \quad \mbox{ and}
$$
$$\|P_{Q,P} \| = \|P_{Q,P_0} \| = (1-\|(1-Q)P_0\|^2)\mrai ,
$$
because $ \ker Q \cap R(P_0) = \{0\}$ \quad \QED

\vglue.3truecm


\begin{thebibliography}{XXXXXX}
\bibitem{[Al1]}  A. Aldroubi, Oblique projections in atomic spaces, Proc. Amer. Math. Soc.
124 (1996), 2051-2060.
\bibitem{[Al2]} A. Aldroubi, Oblique and hierarchical multiwavelet bases, 
Appl. Comput. Harmon. Anal. 4, No.3, 231-263
\bibitem{[AT]}W. N. Anderson and  G. E. Trapp, Shorted operators II, SIAM J. Appl.
Math. 28 (1975), 60-71.
\bibitem{[ACS1]}E. Andruchow,G. Corach and D.  Stojanoff, Geometry of
oblique projections, Studia Math. 137 (1) (1999) 61-79.
\bibitem{[Bo]} R. Bouldin; The product of operators with closed range, Tohoku
Math. J. 25 (1973), 359-363.
\bibitem{[Bu]} D. Buckholtz, Hilbert space idempotents and involutions,
Proc. Amer. Math. Soc. 128 (2000), 1415-1418.
\bibitem {[Ca]} D. Carlson, What are Schur complements, anyway?, Linear
Algebra Appl. 74 (1986), 257-275.
\bibitem {[Co]} R. W. Cottle, Manifestations of the Schur complement, Linear
Algebra Appl. 8 (1974), 189-211.
\bibitem {[Da]} C. Davis, Separation of two linear subspaces, Acta Sci. Math.
(Szeged) 19 (1958), 172-187
\bibitem{[De]} F. Deutsch, The angle between subspaces in Hilbert space,
in "Approximation theory, wavelets and applications" (S. P. Singh, editor), Kluwer,
Netherlands, 1995, 107-130.
\bibitem{[D]}J. Dieudonn\'e, Quasi-hermitian operators, Proc. Internat.
Symp. Linear Spaces, Jerusalem (1961), 115-122.
\bibitem{[Di]} J. Dixmier, Etudes sur les vari\'et\'es et op\'erateurs de Julia, avec
quelques applications, Bull. Soc. Math. France 77 (1949), 11-101.
\bibitem{[Do]}R. G. Douglas, On majorization, factorization and range inclusion
of operators in Hilbert space, Proc. Amer. Math. Soc. 17 (1966) 413-416.
\bibitem{[FW]}P. A. Fillmore and J. P. Williams, On operator ranges,
Advances in Math. 7 (1971) 254-281.
\bibitem{[Ha]} P. R. Halmos, Two subspaces, Trans. Amer. Math.
Soc. 144 (1969), 381-389.
\bibitem{[HN]} S. Hassi,  Nordstr\"om, K.; On projections in a space with an
indefinite metric, Linear Algebra Appl. 208/209 (1994), 401-417.
\bibitem {[Iz]}S. Izumino, The product of operators with closed range and
an extension of the reverse order law, Tohoku Math. J. 34 (1982),
43-52.
\bibitem{[K]} M. G. Krein,  The theory of self-adjoint extensions of semibounded
Hermitian operators and its applications, Mat. Sb. (N. S.) 20 (62)
(1947), 431-495
\bibitem{[L]}P. D. Lax, Symmetrizable linear transformations, Comm. Pure
Appl. Math. 7 (1954), 633-647.
\bibitem{[PW1]}Z. Pasternak-Winiarski, On the dependence of the orthogonal
projector on deformations of the scalar product , Studia Math., 128
(1998), 1-17.
\bibitem{[PEKA]}E. L. Pekarev, Shorts of operators and some extremal
problems, Acta Sci. Math. (Szeged) 56 (1992), 147-163.
\bibitem{[Pt]} V. Ptak, Extremal operators and oblique projections, Casopis
pro pestov\'an\'\i \ Matematiky, 110 (1985), 343-350.
\bibitem{[T]} W. S. Tang, Oblique projections, biorthogonal Riesz bases and
multiwavelets in Hilbert spaces, Proc. Amer. Math. Soc.  128 (2000), 463-473.
\bibitem{[T2]} W. S. Tang,  Oblique multiwavelets in Hilbert spaces, 
Proc. Amer. Math. Soc.,to  appear (electronically published as 
S0002-9939(99)05432-5).
\bibitem{[W1]}H. K. Wimmer, Canonical angles of unitary spaces and perturbations
of direct complements, Linear Alg. Appl. 287 (1999), 373-379.
\bibitem{[W2]}H. K. Wimmer, Lipschitz continuity of oblique projections,
Proc. Amer. Math. Soc. 128 (1999), 873-876.

\end{thebibliography}
\end{document}